\documentclass[aps,prl,twocolumn,showpacs,preprintnumbers,amsmath,amssymb,groupedaddress,superscriptaddress]{revtex4-1}

\usepackage{graphicx}
\usepackage{dcolumn}
\usepackage{bm}
\usepackage{subfigure}

\usepackage{amsmath}
\usepackage{amsthm}
\usepackage{amssymb}

\begin{document}

\affiliation{Center for Applied Mathematics, Cornell University, Ithaca, New York 14853, USA}
\affiliation{School of ORIE, Cornell University, Ithaca, New York 14853 USA}

\date{\today}
\preprint{APS/123-QED}

\title{ Probabilistic Convergence Guarantees for Type II  Pulse Coupled Oscillators}

\author{Joel Nishimura}
\affiliation{Center for Applied Mathematics, Cornell University, Ithaca, New York 14853, USA}
\author{Eric J. Friedman}%
\affiliation{International Computer Science Institute at
Berkeley, California  94704, USA}
\affiliation{School of ORIE, Cornell University, Ithaca, New York 14853 USA}

\begin{abstract}

 We show that a large class of pulse coupled oscillators converge with high probability from random initial conditions on a large class of graphs with time delays. Our analysis combines previous local convergence results, probabilistic network analysis, and a new classification scheme for Type II phase response curves to produce rigorous lower bounds for convergence probabilities based on network density. These bounds are then used to develop a simple, fast and rigorous computational analytic technique. These results suggest new methods for the analysis of pulse coupled oscillators, and provide new insights into the operation of biological Type II phase response curves and also the design of decentralized and minimal clock synchronization schemes in sensor nets.

\end{abstract}


\maketitle

Synchrony in systems of pulse-coupled oscillators (PCOs) is an important feature in physics, biology and engineering.  Synchronization can range from being a pathological breakdown, as in epilepsy \cite{Fisher} to one of vital importance, such as in the proper functioning of the heart's sinoatrial node \cite{Guevara1986,Canavier}, to a framework to understand complex systems \cite{ Strogatz,Timme}.  Additionally, there are attempts to utilize the simplicity of PCO synchronization to synchronize wireless sensor networks \cite{Anna,KeBo,Wang,Pagliari}.  However, many of the idealized models inspired by synchronization are not able to synchronize when the system has a complicated graph structure and time delays -- aspects expected in real physical systems. In order to deal with these issues, previous studies have considered oscillators augmented with memory \cite{KeBo,Hu}, infinite spatial density \cite{Hu} or indegree normalization \cite{Timme, Hu}. While these studies have shown linear stability \cite{Timme}, or other forms of local convergence \cite{KeBo,Hu}, global convergence in these settings has either been shown to be impossible \cite{Timme} or remains unknown.   

 Alternatively, a class of oscillators with Type II phase response curves (PRCs), have been connected to synchronizing behavior theoretically \cite{Aushra,GoelErmentrout} and in nature \cite{Tateno2007,Guevara1986,Anumonwo}. The distinguishing feature of oscillators with Type II PRCs is that an oscillator's phase can either be decreased (inhibited) or increased (excited), depending on the internal state of the oscillator. In this paper we focus on PCOs with a particular class of type II phase response curves, introduced in our previous paper \cite{firstpaper},  which resembles those in nature \cite{Guevara1986,Tateno2007} and are well suited for handling complex topologies and time delay. We also show how leveraging the main theorem from \cite{firstpaper} allows for a computational analytic routine yielding a fast and rigorous estimate of the convergence probability of a system of PCOs.  Furthermore, we provide rigorous lower bounds that guarantee the performance of this computational analytic approach and display how the probability of synchrony converges to $1$ in highly connected graphs. This result is of biological relevance to the situations where synchrony is brought about via Type II PRCs, and is a useful guide for the construction of PCOs in sensor nets.

Previous work found that a class of Type II PRC, denoted ``Stong Type II'' or ``STII'' (described later), could consistently converge to synchrony on fairly complex graphs with time delays \cite{firstpaper}. This convergence was explained by showing that these PCOs would converge to synchrony if their phases were inside a critical range $\rho_0$, essentially showing an $l^\infty$ ball of stability.  This showed that with well-tuned parameters the system is robust to any individual oscillator error or a combination of small errors; explaining the possibility of synchrony, but not the ubiquity of it in numerical simulations.  For example, if  the critical range is $\frac{1}{2}$ of the phase interval, then the probability that a system of $n$ oscillators with uniform random initial conditions starts in the critical regime is ${\frac{1}{2}}^{n}$, which is exponentially small in the system size; however, numerical experiments show that convergence is in fact highly likely and our analysis explains this. 
In particular, we use network analysis to expand on the local understanding of stability, showing that if node indegrees follow a simple scaling law, random initial conditions in any size system are very likely to collapse to the critical range. 

This analysis of STII oscillators sheds light on some of the natural questions regarding Type II phase response curves, and the contribution of excitation and inhibition to synchronization in Type II PRCs.  For example, it's clear that the important aspect of the excitatory end of a Type II phase response curve is that it allows for firing cascades, yet previous analytic results have tended to focus on the importance of inhibition when the system has time delays \cite{firstpaper,Klinglmayr,Timme}. In contrast the result in this paper classifies the excitation in a type II PRC into different discrete classes, each class corresponding to its ability to cascade and a lower bound on the probability convergence.

To understand the strong convergence of STII oscillators, consider the following PCO model:  there are $n$ oscillators on a strongly connected aperiodic directed graph $G$.  Each oscillator $i$'s state is described by phase variable $\phi_i\in [0,1]$, which evolves with natural frequency $d \phi_i/ d t =1$.  When $\phi_i(t)=1$ the oscillator emits a pulse and its phase is reset from $1$ to $0$. This pulse is received by all of $i$'s successors, $S(i)$, time $\tau<0.5$ later.  When a pulse is received, oscillators process this pulse via the phase response curve $f_{ij}(\phi_i)$, where $\phi_j \rightarrow \mathrm{max}(0,\phi_j+f_{ij}(\phi_j))$.  This setup for a PCO can be made to accommodate many of the different types of PCOs examined in the literature. For example, the popular Mirrolo and Strogatz model, which uses a charging curve $V$, can be described by $f_{ij}=V^{-1}(\epsilon +V(\phi_i))-\phi_i$.  

The central goal of this paper is to understand $P_f(G)$, the proportion of phase space that converges to synchrony for a graph $G$ and PRC $f$. In  particular we are interested STII phase response curves defined similarly as in \cite{firstpaper}:
$$   
f(x)=  \left\{
     \begin{array}{lr}
       \le -\mathrm{min}(x,\tau+\kappa)  & : x < B\\
       \ge 0, & : x \ge B
     \end{array}
   \right.
 $$
where $B\in (0,1)$ and $\kappa>0$ are parameters. In \cite{firstpaper} it was shown that if every $f_{ij}$ is a STII PRC, and at some time all oscillator phases are within a critical range $\rho_0< \mathrm{min}(B-\tau,1-B+\tau)$, then the system will converge to synchrony in finite time.  

\begin{figure}
\subfigure{
  \includegraphics[width=\linewidth]{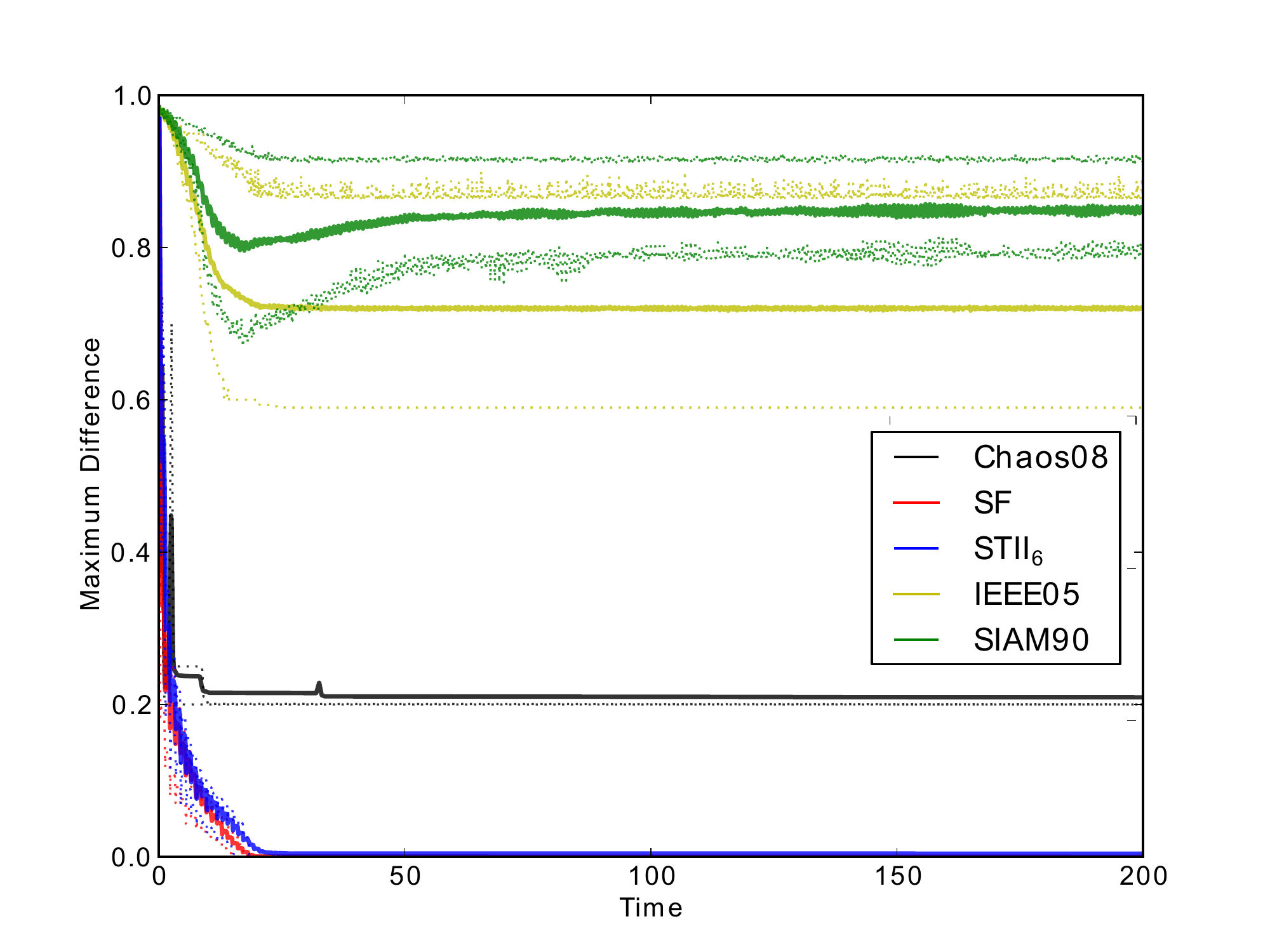}
}
\subfigure{
\includegraphics[width=\linewidth]{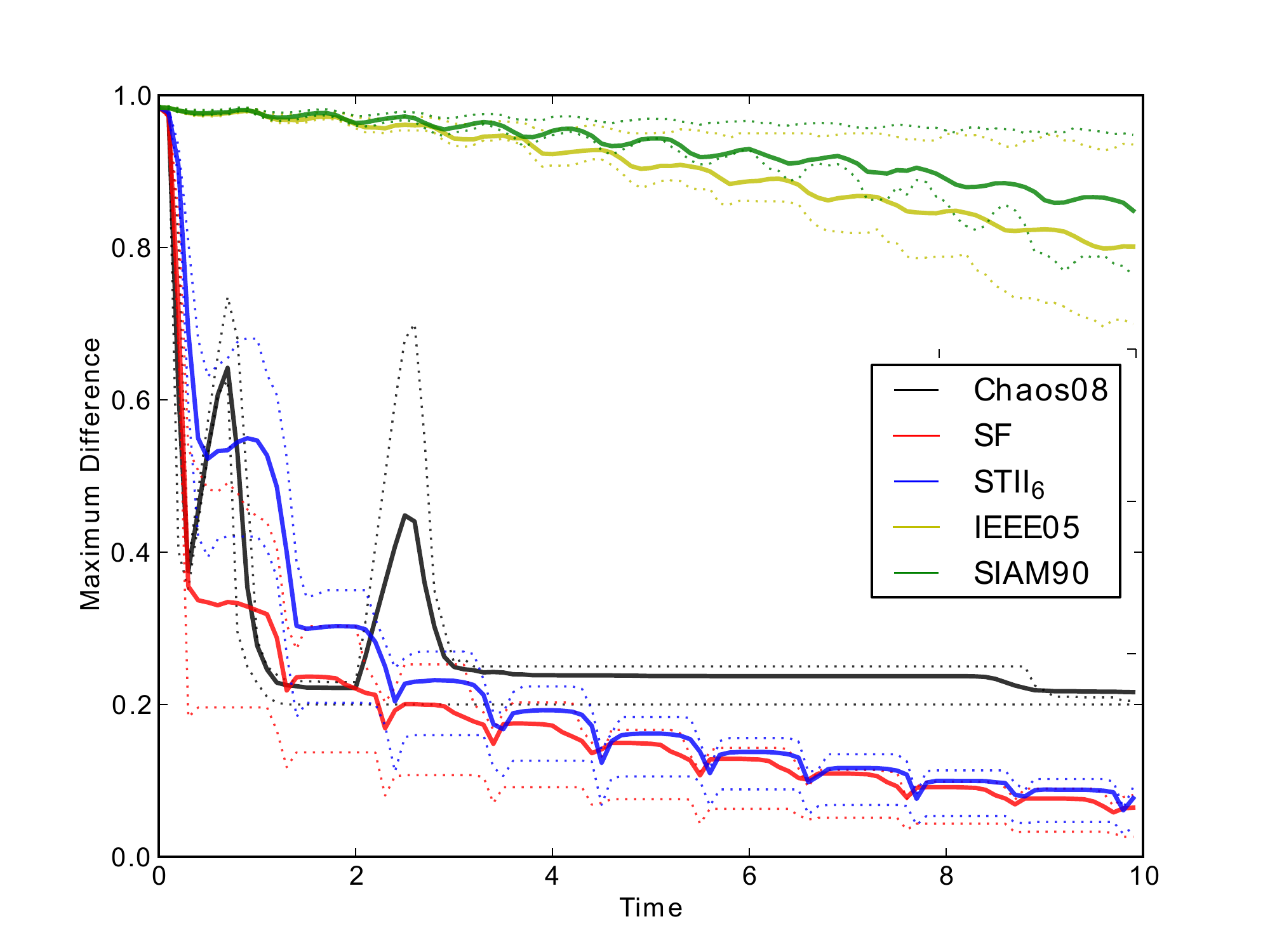}

}
  \caption{Starting from uniform random initial conditions in a system with $400$ nodes, trajectories (mean -- solid line, middle $50\%$ -- between dotted lines) either converge to synchrony or not depending on the PRC.  Notice that two classes discussed in this paper, SF and STII oscillators converge to exact synchrony in finite time while others from Chaos08 \cite{Konishi}, IEEE05 \cite{Anna}, and SIAM90 \cite{Strogatz} do not (this remains true in other measures, not shown). Furthermore, convergence of SF and STII is guaranteed once the maximum difference is smaller than $.5$ by \cite{firstpaper}.}\label{timeConv} 
\end{figure}

To demonstrate the ability of STII oscillators to reach synchrony on complex graphs with time delays, consider Figure \ref{timeConv}, which shows the maximum differences between oscillator phases as a system is integrated for different PRCs and random initial conditions.  Notice that not only are STII curves the only curves that converge, but for most runs, STII curves fall within the critical range in a single time step, despite the fact that the size of the critical range is exponentially small in probability space.  As will be shown, a large portion of the basin of attraction of synchrony in STII oscillators can be described by this rapid convergence to the critical range.  Furthermore, this convergence to the critical range arises from a fundamentally different mechanism, and relies on different properties of the PRC than the convergence inside the critical range.

To understand this basin analytically, consider first, the most extreme STII PRC, the ``strong firing'' (SF) PRC where $f^{SF}(x)=-x$ for $x<B$ and  $f^{SF}(x)=1-x$ otherwise.  Notice, the response $f^{SF}$ gives to any signal causes an oscillator to have phase $0$, where signals received after $B$ also cause it to fire.  This brief characterization of the SF PRC allows for a quick analytic description of one way in which the SF PRC converges rapidly to the critical range.

The key insight is that if for every SF oscillator $i$, $i$ receives a signal or fires in a small window of time, denoted as event $E_i$, then every oscillator will be reset to phase $0$, and thus all phases will be within the critical range.  We will show that the window can be as large as $1-s$, where $s=\mathrm{max}(B,1-B+2\tau)$. Notice, that if every oscillator receives a signal or fires at some time in $[\tau,1-s+\tau]$, then for all $i$, $\phi_i (1-s+\tau )\le 1-s \le B-2\tau$.  Since only oscillators with phases greater than $B$ can generate signals, no new signals are being sent.  Thus by time $\tau$ later all signals will have arrived at their destinations, giving that for all $i$, $\phi_i (1-s+2\tau) \le 1-s+\tau \le \rho_0$ which puts all oscillators in the critical range with no signals enroute.  

Therefore, the probability that the SF system converges can be bounded by the likelihood that every oscillator receives a signal in a window of size $1-s$ time, $P(\cap_i ^n E_i)$, which can then be bounded by using node degrees, $d_i$.  If at time $0$, $\phi_j(0)\in [s,1]$ then all the successors of $j$ will necessarily receive a signal in time $[\tau,1-s+\tau]$ (since $s\ge B$, $j$ is in the excitatory regime) and if $\phi_j(0)\in [s-\tau,1-\tau]$ then $j$ must fire in  $[\tau,1-s+\tau]$.  (Note: that this is now very similar to the probability that a random subset of the graph will dominate it, connecting it to some sensor net protocols used to find a Connected Dominating Set \cite{Ibanez,Dai} and the study of dominating sets in general \cite{Chellali}). Thus, for uniform random initial conditions a simple bound on the probability of any oscillator $i$ receiving a signal or firing in $[\tau,1-s+\tau]$, is simply the complement of all $i$'s predecessors having phases in $[0,s)$ and $\phi_i \in [0,s-\tau] \cup [1-\tau,1]$; yielding $P(E_i)\ge 1- s^{d_i+1}$.  

Notice, that the probability of a node $i$ failing to receive a signal in the $1-s$ time window is exponentially small in that node's indegree.  These probabilities can be aggregated using the Union Bound, giving that: $P(\cup_i ^n E_i ^c )\le \Sigma_i ^n s^{d_i+1}$ and thus $P(\cap_i ^n E_i)\ge 1- \Sigma_i ^n s^{d_i+1}$. Alternatively, a slightly stronger bound can be found using the fact that each $E_i$ is positively correlated, or that the number of nodes dominated is a submodular function of node subsets, giving the probability of convergence,
$P_{SF}(G)\ge P(\cap_i ^n E_i )\ge \Pi_i ^n (1-s^{d_i+1}) .$ Thus the deterministic bound from \cite{firstpaper} has been used to create a statement about convergence from random initial conditions.  


%

This result immediately gives a number of interesting corollaries.  For example, let $\delta_n(p)$ be the minimum indegree such that the $P_{SF}(G)>p$ then $\delta_n(p)\ \le \ln (1-p)/\ln (s) - \ln (n)/\ln (s) ,$ which is logarithmic in the system size.  To give a sense of the constants, the minimum value for $s$ occurs when $B=s=.5+\tau$, and thus to ensure a $95\%$ convergence rate in a systems with a time delay $5\%$ of the period: $\delta_n(0.95)\le 5.02+1.68 \ln (n)$; a result that holds for any $n$.

This result can also be used to make statements about the convergence of SF oscillators on random graphs. Take an Erd\H{o}s-R\'{e}nyi  random graph, $G(n,\hat{p})$, where edges are created with independent probability $\hat{p}$.  In this case an application of Chernoff's inequality shows that if $\hat{p}\ge \frac{\ln (n)}{n} g(s,\gamma)$ for a function $g$ of $s$ and some number $\gamma>0$ then as $n \rightarrow \infty$,  the probability of synchrony $P_{SF}(G) \to (1-1/n)e^{1/n^{1-\gamma}} \to 1$. 
 Notice, that this requirement on $\hat{p}$ is only a constant multiple of that required for $G(n,\hat{p})$ to be connected, which  asymptotically, occurs when $\hat{p}\ge (1+\epsilon)\frac{\ln{n}}{n}$. Thus, the degree requirements grow reasonably with $n$, and furthermore, since convergence requires connectedness, our rigorous bound is a constant factor approximation of the actual required degree.

Similarly, one can show asymptotic bounds on random geometric graphs, constructed by positioning nodes uniformly at random on the unit $\bar{d}$ dimensional torus and connecting any nodes within some radius $r$. If $r$ is chosen so that the expected degree $r^{\bar{d}} n\theta=c \ln(n)$, (where $\theta$ is the volume of a $\bar{d}$ dimensional unit ball) then utilizing results describing the minimum degree in random geometric graphs, \cite{penrose} shows that the system will converge to synchrony as $n \to \infty$ so long as $c$ is the greater of the solutions to: $\frac{1}{c}=1+\frac{2}{c\ln{s}}-\frac{2}{c\ln{s}}\ln( \frac{-2}{c\ln{s}})$.  Again, a constant factor of logarithmic growth in the expected degree gives convergence guarantees.  

\begin{figure}
\includegraphics[width=\linewidth]{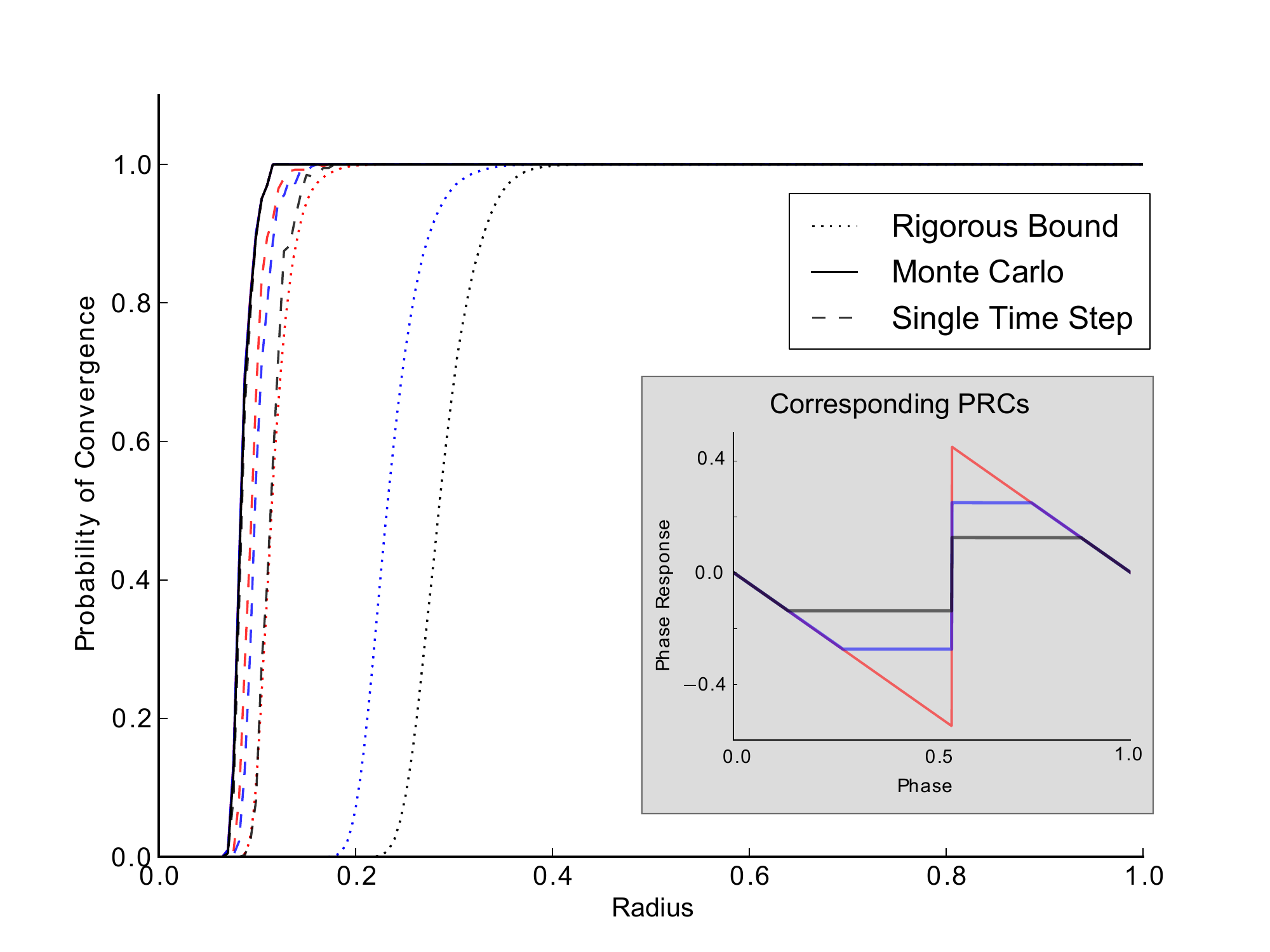}\\

  \caption{ Probability of convergence  for a $400$ node random geometric graph as a function of radius for SF (red), $STII_{4,0}$ (blue) and a $STII_{7,0}$ (black) PCOs.  Numerical results (solid) suggest that  all three oscillators systems transition to synchrony at the same value of $r$. Dotted lines show an analytic lower bound and dashed lines show the numerical single time step bound. 
  }\label{convG} 
\end{figure}

Thus far, we have shown analytically that sufficiently dense systems of SF PCOs will, with high probability, converge in a single time step to the critical range, and consequently, will converge to exact synchrony in finite time. However, synchrony is a real phenomenon, and the constant in our bound  may be important in certain applications. Fortunately, the analysis of the lower bound can be extended in a simple ``computational-analytic'' manner, providing a rigorous computational assisted bound.  Figure \ref{convG} shows the relationship between indegree,  the analytic lower bound, the numerical bound and  two computational analytic bounds.  
The first, simulates the system for {\em one period} and then applies the deterministic convergence result.  This is fast, works surprising well and  is far more computationally efficient than fully integrating the system to convergence (a time that would scale with a graph's aperiodic diameter). For situations where precise estimates of $P_{SF}(G)$ are important one uses the computational analytic approach by running multiple single time step Monte Carlo trials and checking if the phases fall within the critical range. In a graph with $m$ edges such a routine can be implemented by an event based simulation, and thus can run in $O(m\log{m})$ time. Whereas typically integration time scales with system size, our analytic bound guarantees that this computational analytic routine remains viable as the system size increases.

\begin{figure}
  \includegraphics[width= \linewidth]{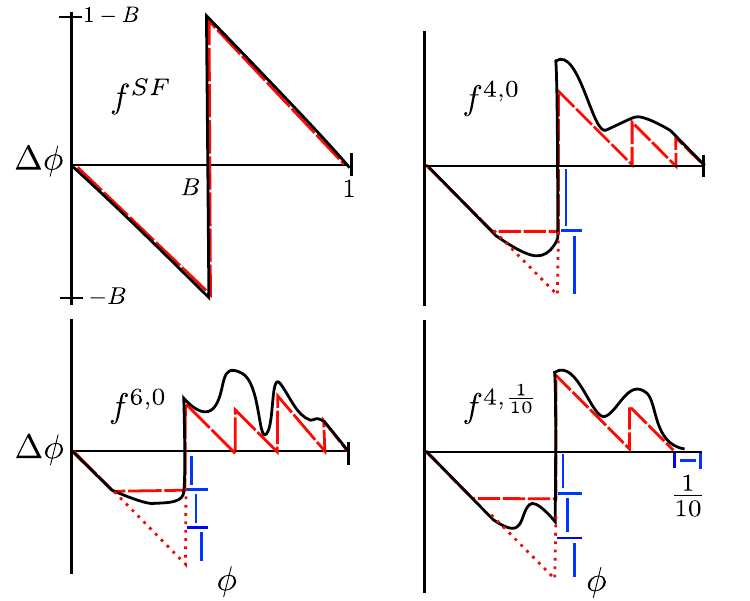}\\
  \caption{ Placing a sawtooth function underneath the excitatory portion of an $STII$ PRC provides an upper bound on the number of excitations that will cause the curve to fire, whereas the scale of a curve's inhibition is measured in proportion to $B$. }\label{saw} 
\end{figure}

We now consider more general PRCs, showing that many STII PRCs will also have $\delta_n(p)=O(\mathrm{log}(n))$ and thus will also have corresponding guarantees for random graph models and computational analytic routines. The key feature of the SF PRC was that a single signal causes an oscillator to reset or fire. The arguments made for an SF oscillator can be modified to allow for oscillators that require multiple signals to reset or fire.  Consider the sub class $STII_{k,\eta}$ which, as opposed to requiring $1$ signal, will require receiving at least $k$ signals within $1-s-\eta$ time to reset or fire.  Figure \ref{saw}, displays several different PRCs from $STII_{k,\eta}$ for different $k$ and $\eta$.

If for independent initial conditions, $P(\phi_i(0) \in [0,s+\eta)) = q_i$, then the probability that an oscillator $i$ has less than $k$ neighbors ready to fire is the sum of the binomial distribution $B(d_i,q_i)$ from $0$ to $k-1$. Using a well known bound based on Hoeffding's inequality yields that the probability that $i$ has at less than $k$ neighbors ready to fire is less than $e^{\frac{-(q_id_i-k-1)^2}{2q_id_i} }$.  Taking the complement and using the same results of independence gives that the probability of the system synchronizing is, $$ P_{f_{k,\eta}}(G) \ge \Pi_{i=1}^n (1-e^{\frac{-(q_id_i-k-1)^2}{2q_i d_i} } ).$$

Alternatively, let $c_n=\ln(n)-\ln(1-p)$, then the system will converge with at least probability $p$ if  for all $i$, the expected number of firing neighbors, $d_i q_i \ge k -1+ c_n +\sqrt{c_n^2+4(k-1)c_n}$.  Since, for fixed $k$ this result also scales $O(\ln(n))$ then the random graph results in the SF case have analogs of the same order: $\hat{p} = O(\frac{\ln(n)}{n})$ for Erd\H{o}s-R\'{e}nyi  and $r=O( \frac{\ln (n)}{n^{1/\bar{d}}})$ for random geometric graphs.

Determining if a phase response curve is a member of $STII_{k,\eta}$ involves two steps: first, classifying the strength of the inhibitory section, and second, the strength of the excitation. 

We say that an oscillator $i$ is $h$-inhibitory if receiving $h$ signals in the inhibitory region over some span of time $[t_0,t_0+s']$, $s'=1-s-\eta$ forces $\phi_i(t_0+s')<s'$.  For example, a sufficient condition for a PRC $f$ to be $h$-inhibitory is if $f(x)< -\mathrm{min}(\frac{B}{h},x)$, for $x<B$.  For such an $f$, $h$ signals causes the oscillator either to be reset to $0$, or to be inhibited by at least $B$, giving that $\phi(t_0+s')<s'$.  

Similarly, an oscillator $i$ is $(k-h+1)$-excitatory if receiving $k-h+1$ signals in the excitatory region, in some time $[t_0,t_0+s']$, forces an oscillator to fire before $t_0+s'+\eta$.  As seen in figure \ref{saw}, a sufficient condition for $(k-h+1)$-excitability is that the PRC is greater than a saw tooth with slope $-1$ when $x\in[B,1-\eta]$.  

If an oscillator is both $h$-inhibitory and $(k-h+1)$-excitatory then it is a member of $STII_{k,\eta}$.  These results can also be used to guarantee the performance of a similar computational analytic routine. The performance of a $STII_{7,0}$ and a $STII_{4,0}$ as well as their analytic guarantees can be seen in figure \ref{convG}.

Finally, it is worth noting that in many systems, such as systems of neurons, edges are often weighted and the impact that different neighbors have varies drastically \cite{Timme}.  The results in this paper also extend to a weighted version, where each edge has weight $w_{ij}$ and weights are interpreted by the formula: $\hat{f}_{ij}(x)=\mathrm{max}(-w_{ij},f_{ij})$ for $x<B$ and $f_{ij}(x)=\mathrm{min}(w_{ij},f_{ij})$ for $x>B$, where $w_{ij}$ acts as a constraint of the phase response curve.  The above formula for the $P_f(G)$ remains true so long as for each $i$, $\sum_j w_{j,i} \ge \tau$ and for each node $i$ there are $d_i$ nodes $j$ such that $f_{ij} \in F_{h,k}$.  

In such a case, if $k$ increases as $O(\ln(n))$ or less then so does $\delta_n(p)$.  Furthermore, if $k\rightarrow\infty$ then the requirements on the phase response curve shrink to simply requiring that $f'(0)=-1$, and $f(x)<0$ for $x<B$ and $f(x)>0$ for $x>B$ and that $f$ is continuous everywhere except $B$ where $\lim_{x  \to B^-}<-\epsilon$ and $\lim_{x \to B^+}>\epsilon$. Thus for very large systems our results show convergence for a very general class of type II oscillators.  However, when comparing to results for ``weakly coupled oscillators'' one should recall the slightly different requirement, that in those cases, $\sum_i w_{ij} =\epsilon$ and the weights are multiplicative: $f_{ij}(\phi)=w_{ij}f(\phi)$.

In summary, we have shown how the local convergence of STII pulse coupled oscillators to synchrony can be extended probabilistically, relating graph density and phase response curve structure to a rigorous lower bound on the probability of convergence. Applying this lower bound to random graph models shows that the expected node degree beyond which synchronization is very likely is a constant multiple of the percolation threshold. Therefore a computational analytic scheme that simply sampled single time steps is a constant factor approximation to a sampling routine that integrated for infinite time. An extension with edge weights was also discussed.

%
%
%
%
%

\bibliography{OscReferences2}

\end{document}